\documentclass[11pt]{amsart}
\usepackage{amssymb}
\newfont{\cyr}{wncyr9}
\newtheorem{theorem}{Theorem}[section]
\newtheorem{corollary}[theorem]{Corollary}
\newtheorem{lemma}[theorem]{Lemma}
\newtheorem{question}[theorem]{Question}

\newtheorem{historical background}[theorem]{Historical Background}
\newtheorem{outline}[theorem]{Outline}
\newtheorem{remarks}[theorem]{Remarks}
\newtheorem{remark}[theorem]{Remark}

\newtheorem{notation}[theorem]{Notation}

\def\cc {{\mathfrak c}}
\def\ttt {{\mathfrak t}}
\def\uu {{\mathfrak u}}

\def\AA {{\mathbb A}}

\def\PP {{\Bbb P}}

\def\RR {{\Bbb R}}
\def\TT {{\Bbb T}}

\def\ZZ {{\Bbb Z}}
\def\sA {{\mathcal A}}
\def\sB {{\mathcal B}}
\def\sC {{\mathcal C}}

\def\sH {{\mathcal H}}

\def\sM {{\mathcal M}}

\def\sP {{\mathcal P}}

\def\sS {{\mathcal S}}
\def\sT {{\mathcal T}}

\def\cf {\mathrm{cf}}

\def\Fac {\mathrm{Fac}}
\def\Hom {\mathrm{Hom}}

\def\graph {\mathrm{graph}}

\def\mod {\mathrm{mod}}

\def\< {{\langle}}
\def\> {{\rangle}}

\begin{document}
\title{Making group topologies with, and without, convergent sequences}
\author{W. W. Comfort}
\address{W. W. Comfort: Department of Mathematics, 
Wesleyan University, Middletown, CT 06459}
\email{wcomfort@wesleyan.edu}
\author{S. U. Raczkowski}
\address{S. U. Raczkowski and F. J. Trigos-Arrieta: 
Department of Mathematics, California State University, 
Bakersfield, Bakersfield, CA, 93311-1099}
\email{racz@csub.edu}
\author{F. J. Trigos-Arrieta}
\email{jtrigos@csub.edu}
\begin{abstract}
(1) Every infinite, Abelian compact (Hausdorff) group $K$ admits 
$2^{|K|}$-many dense, non-Haar-measurable subgroups of cardinality $|K|$. 
When $K$ is nonmetrizable, these may be chosen to be pseudocompact. 

(2) Every infinite Abelian group $G$ admits a family $\mathcal{A}$ of 
$2^{2^{|G|}}$-many pairwise nonhomeomorphic totally bounded group  
topologies such that no nontrivial sequence in $G$ converges in any of 
the topologies $\mathcal{T}\in\mathcal{A}$. (For some $G$ one may 
arrange $w(G,\mathcal{T})<2^{|G|}$ 
for some $\sT\in\mathcal{A}$.)

(3) Every infinite Abelian group $G$ admits a family $\mathcal{B}$ of 
$2^{2^{|G|}}$-many pairwise nonhomeomorphic totally bounded group 
topologies, with\break $w(G,\mathcal{T})=2^{|G|}$ for all 
$\mathcal{T}\in\mathcal{B}$, such  that some fixed faithfully indexed 
sequence in $G$ converges to $0_G$ in each $\mathcal{T}\in\mathcal{B}$.
\end{abstract}
\thanks{Portions of this paper were presented by the first-listed author 
at the 2004 Annual Meeting of the American Mathematical Society (Phoenix, 
January, 2004).}
\thanks{The second listed author acknowledges partial support from 
the University Research Council at CSU Bakersfield. She also wishes to 
thank Mrs.~Mary Connie Comfort for her encouragement, without which 
this paper would never see the daylight. Thank you.}
\maketitle

\section {Introduction}

\begin{historical background}
{\rm Not long after M. H. Stone and E. \v{C}ech associated with each
Tychonoff space $X$ its maximal compactification $\beta(X)$ (the
so-called Stone-\v{C}ech compactification), it was noted, denoting by
$\omega$ the countably infinite discrete space, that $\beta(\omega)$
contains no nontrivial convergent sequence. This observation stimulated
Efimov~\cite{efimov69} to pose in 1969 a question which in its full
generality remains unsolved today: Does every compact, Hausdorff space
contain either a copy of $\beta(\omega)$ or a nontrivial convergent
sequence? (In models of $\diamondsuit$ the answer is
negative~\cite{fedorcuk76}.
See \cite{shakh92} for several additional relevant references.) The
present paper is concerned with topological groups. In that context, a
correct and natural companion to Efimov's question is this: 
Given a class $\sC$ of
topological groups, does every group in $\sC$ contain a
nontrivial convergent sequence? There is an extensive literature on
questions of this form. Here are some samples
of both positive and negative results.

\underline{Positive} (a)~According to a result to which
{\v{S}}apirovski{\u{\i}}, Gerlits and Efimov have contributed (see
\cite{shak94} for historical details and for an ``elementary" proof), every
infinite compact group $K$ contains topologically a copy of the
generalized Cantor set $\{0,1\}^{w(K)}$, hence contains a convergent sequence;
(b)~Assuming GCH, Malykhin and Shapiro~\cite{malyshap} showed that every totally
bounded group $G$ with $w(G)<(w(G))^\omega$ contains a nontrivial
convergent sequence; (c)~Raczkowski~\cite{raczphd},~\cite{racz:02}
and others \cite{Barb} have shown that
for every suitably fast-growing sequence $x_n\in\ZZ$ there is a totally
bounded group topology on $\ZZ$ with respect to which $x_n\rightarrow0$.

\underline{Negative}. Glicksberg~\cite{glickii} showed that when a locally
compact Abelian group $(G,\sT)$ is given its associated Bohr topology
(that is, the weak topology induced on $G$ by $\widehat{(G,\sT)}$), no
new compact sets are created; in particular, as shown earlier by
Leptin~\cite{leptin55a}, the topology induced on a
(discrete) Abelian group by $\Hom(G,\TT)$ has no infinite compact
subsets, in particular has no nontrivial convergent sequence; (b)~there
are infinite pseudocompact topological groups containing no nontrivial
convergent sequence~\cite{sirota}. 

Perhaps the most celebrated unsolved question in this area of
mathematics is this: Does there exist in ZFC a countably compact topological group
with no nontrivial convergent sequence? (Many examples are known in
augmented axiom systems. See for example the constructions of van
Douwen~\cite{vdi}, of Hart and van Mill~\cite{hartvmii}, and of
Tomita~\cite{tomita96},~\cite{tomita99}, and see also \cite{diksha03} for a
characterization, in a forcing model of ZFC~+~CH with $2^\cc$
``arbitrarily large", of those Abelian groups which admit a hereditarily
separable pseudocompact (alternatively, countably compact) group topology with no infinite
compact subsets.)
}
\end{historical background}

\begin{outline}
{\rm The present work achieves results which in a certain direction
may legitimately be called
optimal: We show that every infinite Abelian group
$G$ admits the maximal number, namely
$2^{2^{|G|}}$, of totally bounded group topologies with no nontrivial convergent
sequence (Theorem~4.1); and also, the same number of totally bounded group topologies in each of
which some nontrivial sequence (fixed, and chosen in advance) does
converge (Theorem~5.5). An elementary cardinality argument shows that
the various topological groups $(G,\sT)$ may be chosen to be pairwise nonhomeomorphic
as topological spaces.
}
\end{outline}

\begin{notation} 
{\rm
The symbols $\kappa$ and $\alpha$ denote infinite
cardinals; $\omega$ is the least infinite cardinal, and $\cc:=2^\omega$.
For $S$ a set we write $[S]^\kappa:=\{A\subseteq S:|A|=\kappa\}$; the symbols
$[S]^{<\kappa}$ and $[S]^{\leq\kappa}$ are defined analogously.
$\ZZ$ and
$\RR$ denote respectively the group of integers and the group of
real numbers, often with their usual (metrizable) topologies, and
$\TT:=\RR/\ZZ$. General groups $G$ are written multiplicatively,
with identity $1=1_G$,
but groups $G$ known or hypothesized to be Abelian are written
additively, with identity $0=0_G$. We identify each finite cyclic group
with its copy in $\TT$; in particular for $0<n<\omega$ we write
$\ZZ(n):=\{\frac{k}{n}:0\leq k<n\}\subseteq\TT$. We denote by
$\PP$ the set of primes.

We work exclusively with Tychonoff spaces, i.e., with completely
regular, Hausdorff spaces. A topological group which is a Hausdorff space is
necessarily a Tychonoff space~\cite{ehkr:63}(8.4).
A topological group $(G,\sT)$ is said to be {\it totally bounded} if for
$\emptyset\neq U\in\sT$ there is $F\in[G]^{<\omega}$ such that $G=FU$.
It is a theorem of Weil~\cite{Weil:37} that a topological group is totally bounded if and only if $G$
embeds densely into a compact group; this latter group, unique in an
obvious sense, is called the {\it Weil completion} of $G$ and is denoted
$\overline{G}$. 

Given a topological group $G=(G,\sT)$, the symbol $\widehat{G}=\widehat{(G,\sT)}$
denotes the set of continuous homomorphisms from $G$ to the (compact)
group $\TT$. A group $G$ with its discrete topology is written $G_d$, so
$\widehat{G_d}=\Hom(G,\TT)$; this is a closed subgroup of the compact group
$\TT^G$. It is easily seen, as in \cite{wckr:64}(1.9), that for a subgroup $H$ of
$\Hom(G,\TT)$ these conditions are equivalent: (a)~$H$ separates points
of $G$; (b)~$H$ is dense in the compact group $\Hom(G,\TT)$.

For $G$ Abelian we denote by $\sS(G)$ the set of point-separating subgroups
of $\Hom(G,\TT)$, and by $\ttt(G)$ the set of (Hausdorff) totally
bounded group topologies on
$G$. Theorem~1.5 below describes, for each infinite Abelian group $G$,
a useful order-preserving bijection between $\sS(G)$ and $\ttt(G)$.

The rank, the torsion-free rank, and (for $p\in\PP$) the $p$-rank of an
Abelian group $G$ are denoted $r(G)$, $r_0(G)$, and $r_p(G)$,
respectively. We have $r(G)=r_0(G)+\sum_{p\in\PP}\,r_p(G)$, and when
$|G|>\omega$ we have $|G|=r(G)$ (cf. \cite{fuch:70}(\S16)
or \cite{ehkr:63}(Appendix~A)). Following those sources we write
$$G_p:=\bigcup_{k<\omega}\,\{x\in G:p^k\cdot x=0\}.$$

We write $X=_h Y$ if $X$ and $Y$ are homeomorphic topological spaces,
and we write $G\simeq H$ if $G$ and $H$ are isomorphic groups; it is
important to note that $X=_h Y$ conveys no information about the
algebraic structure of $X$ or $Y$ (if any), and $G\simeq H$ conveys no
information about the topological structure of $G$ or $H$ (if any).

We assume familiarity on the reader's part with the essentials of Haar
measure $\lambda=\lambda_K$ on a locally compact group $K$. The set of
Borel sets, and the set of $\lambda$-measurable sets, are denoted
$\sB(K)$ and $\sM(K)$, respecively. Of course $\sB(K)\subseteq\sM(K)$.
Our
convention is that Haar measure is {\it complete} in the sense that if
$S\in\sM(K)$ with $\lambda(S)=0$, then each
$A\subseteq S$ satisfies $A\in\sM(K)$ (with $\lambda(A)=0$). We assume
for simplicity that if $K$ is compact then $\lambda$ is {\it normalized}
in the sense that $\lambda(K)=1$.

The following three theorems are essential to our argument. Therem~1.4(a)
is due to Steinhaus~\cite{steinhaus20} when $G=\RR$ and to
Weil~\cite{weilii}(p.~50) in the general case; see Stromberg~\cite{strom72} for a
pleasing, efficient proof. Parts (b) and (c) are the immediate
consequences which we use here frequently.
}
\end{notation}

\begin{theorem} 
\cite{steinhaus20}, \cite{weilii}. Let $K$ be a locally compact
group and let $S$ be a $\lambda$-measurable subset of $K$ with
$\lambda(S)>0$. Then
\begin{itemize}
\item[(a)] the difference set $SS^{-1}$ contains a neighborhood of $1_K$;
\item[(b)] if $S$ is a subgroup of $K$ then $S$ is open and closed in 
$K$; and
\item[(c)] if $S$ is a dense subgroup of $K$ then $S=K$.~$\square$
\end{itemize}
\end{theorem}

\begin{theorem} \cite{wckr:64}. Let $G$ be an Abelian group.
\begin{itemize}
\item[(a)] For every $H\in\sS(G)$, the topology $\sT_H$ induced on $G$ by $H$
is a (Hausdorff) totally bounded group topology such that $w(G,\sT_H)=|H|$;
\item[(b)] if $(G,\sT)$ is totally bounded then $\sT=\sT_H$ with
$H:=\widehat{(G,\sT)}\in\sS(G)$.~$\square$
\end{itemize}
\end{theorem}

It is clear with $G$ as in Theorem 1.5 that distinct $H_0,H_1\in\sS(G)$
induce distinct topologies $\sT_{H_0}, \sT_{H_1}\in\ttt(G)$---indeed
$\widehat{(G,\sT_{H_0})}=H_0\neq H_1=\widehat{(G,\sT_{H_1}})$; thus the bijection
$\ttt(G)\leftrightarrow\sS(G)$ given by $\sT_H\leftrightarrow H$ is
indeed order-preserving.

\begin{theorem}  
\cite{CTWu:93}. Let $G$ be an infinite Abelian group
with $K:=\Hom(G,\TT)=\widehat{G_d}$, let
$(x_n)_n$ be a faithfully indexed sequence in $G$, and let
$$A:=\{h\in K:h(x_n)\rightarrow 0\}.$$
Then $A\in\sM(K)$, and $\lambda(A)=0$.
\end{theorem}

\begin{proof}
${\rm [Outline]}$. According to the duality theorem of
Pontrjagin~\cite{PvKd:34} and van~Kampen~\cite{VK:35}, the map
$G_d\twoheadrightarrow\widehat{\widehat{G_d}}=\widehat{\Hom(G,\TT)}$
given by $x\rightarrow\widehat{x}$ (with $\widehat{x}(h)=h(x)$ for $x\in
G$, $h\in K$) is a bijection. Writing
$$A_{n,m}:=\{h\in K:|\widehat{x_n}(h)-1|<\frac{1}{m}\},$$
the relation
$A=\bigcap_{m<\omega}\,\bigcup_{N\geq m}\,\bigcap_{n\geq N}\,A_{n,m}$
expresses $A$ as a $G_{\delta\sigma\delta}$-subset of the
compact group $K$, so $A\in\sB(K)\subseteq\sM(K)$. If $G$ is
torsion-free, a condition equivalent to the condition that $\Hom(G,\TT)$
is connected (cf.~\cite{ehkr:63}(24.25)), then a
reference to Theorem~1.4 completes
the proof: the condition $\lambda(A)>0$ would imply $A=\Hom(G,\TT)$,
so that $x_n\rightarrow0$ in the Bohr topology of $G_d$, contrary to
the theorem of Leptin and Glicksberg cited above. We
refer the reader to \cite{CTWu:93} for the proof (for general Abelian
$G$) that $\lambda(A)=0$.
\end{proof}

In what follows we will frequently invoke this simple algebraic fact.

\begin{theorem}   
Let $K$ be an Abelian group, let $H$ be a subgroup of
$K$ of index $\alpha>\omega$, and let
$$\sH:=\{S:H\subseteq S\subseteq K,\ \text{$S$ is a proper subgroup of 
$K$}, |S|=|K|\}.$$
Then $|\sH|=2^\alpha$.
\end{theorem}

\begin{proof}
The inequality $\leq$ is obvious. We have
$|K/H|=r(K/H)=\alpha>\omega$, so algebraically
$K/H\supseteq\bigoplus_{\xi<\alpha}\,C_\xi$ with each $C_\xi$ cyclic. Let
$\phi:K\twoheadrightarrow K/H$ be the canonical homomorphism, and for
$A\in[\alpha]^\alpha\backslash\{\alpha\}$ set
$H_A:=\phi^{-1}(\bigoplus_{\xi\in A}\,C_\xi)$. The map
$[\alpha]^\alpha\backslash\{\alpha\}\rightarrow\sH$ given by
$A\rightarrow H_A$ is an injection, so $|\sH|\geq2^\alpha=2^{|K/H|}$, as
asserted.
\end{proof}

Theorem 1.6 explains our interest in the existence of (many)
point-separating nonmeasurable subgroups of compact Abelian groups $K$ (of
the form $K=\Hom(G,\TT)$): each such $H\in\sS(G)$ will induce on $G$ a
totally bounded group topology without nontrivial convergent sequences. In order to show that
such $K$ admit $2^{|K|}$-many such subgroups, we find it convenient to treat
separately the metrizable case (that is, $w(K)=\omega$) and the
nonmetrizable case ($w(K)>\omega$). We do this in Sections 2 and 3
respectively. 

\section{Many Nonmeasurable Subgroups: The Metrizable Case}

For simplicity, and because it suffices for our applications, we take
the groups $K$ and $M$ in Lemma~2.4 and Theorem 2.2 to be compact;
the reader may notice that this hypothesis can be significantly relaxed.
Indeed both groups are Abelian and $M$ is metrizable
in our applications, but
since those hypotheses save no labor we omit them for now. 

\begin{lemma}  
 Let $K$ and $M$ be compact groups with Haar measures $\lambda$
and $\mu$ resectively, let $\phi:K\twoheadrightarrow M$ be a continuous
surjective homomorphism, and define $\Phi:\sB(M)\rightarrow\sB(K)$ by
$\Phi(E):=\phi^{-1}[E]$. Then the map
$m:=\lambda\circ\Phi:\sB(M)\rightarrow[0,1]$ satisfies $m=\mu|\sB(M)$.
\end{lemma}

\begin{proof} According to \cite{ehkr:63}(15.8), and using the numbering system
there, it is enough to show that 
\begin{itemize}
\item[(iv)] 
$m(C)<\infty$ for compact $C\in\sB(M)$; (v)~$m(U)>0$ for some open 
$U\in\sB(M)$;
\item[(vi)]
$m(a+F)=m(F)$ for all $a\in M$, $F\in\sB(M)$; and
\item[(vii)] 
$m(U)=\sup\{m(F):F\subseteq U, F$ is compact$\}$ for open $U\subseteq M$, 
and $m(E)=\inf\{m(U):E\subseteq U, U$ is open$\}$ for
$E\in\sB(M)$.
\end{itemize}
The verifications are routine and will not be reproduced here.
In addition to \cite{ehkr:63},
the reader seeking hints might consult \cite{Ha}(63C and
64H), or \cite{Ha}(52G and 52H).
\end{proof}

\begin{theorem} 
Let $K$ and $M$ be compact groups with Haar measures
$\lambda$ and $\mu$ respectively, and let $\phi:K\twoheadrightarrow M$
be a continuous, surjective homomorphism. If $D$ is a dense,
non-$\mu$-measurable subgroup of $M$, then $H:=\phi^{-1}(D)$ is a dense,
non-$\lambda$-measurable subgroup of $K$.
\end{theorem}

\begin{proof} $\phi$ is an open map \cite{ehkr:63}(5.29), so $H$ is dense in $K$.
Suppose now that $H\in\sM(K)$, so that either
$\lambda(H)>0$ or $\lambda(H)=0$. If $\lambda(H)>0$ then
$H=K$ by Theorem~1.4(c) so $D=\phi[K]=M$, a
contradiction. If ($H\in\sM(K)$ and) $\lambda(H)=0$ then since
$\lambda$ is (inner-) regular there is a sequence $K_n $ ($n<\omega$) of
compact subsets of $K\backslash H$ such that
$\lambda(\bigcup_n\,K_n)=\lambda(K\backslash H)=1$. We write
$M_n:=\phi[K_n]$ and $\widetilde{K_n}:=\phi^{-1}(M_n)$. Then
$K_n\subseteq\widetilde{K_n}\subseteq K\backslash H$ and from Lemma~2.1
we have
$$\mu(\bigcup_n\,M_n)
= \lambda(\bigcup_n\,\widetilde{K_n})
\geq \lambda(\bigcup_n\,K_n)
=1,$$
so $\mu(\bigcup_n\,M_n)=1$ and hence $D\in\sM(M)$ with
$\mu(D)=0$, a contradiction.
\end{proof}

Our goal is to show that every (infinite) compact Abelian metrizable
group contains a dense, nonmeasurable subgroup of index $\cc$. We treat some
special cases first. In what follows we denote the torsion subgroup of an
Abelian group $K$ by $t(K)$, for $s\in K$ and $0\neq n\in\ZZ$ we write
$[\frac{s}{n}]=\{x\in K:nx=s\}$, and for a subgroup $S$ of $K$ we set
$$\mbox{div}(S)
:= \bigcup\{[\frac{s}{n}]:s\in S, 0\neq n\in\ZZ\}.$$
When $[\frac{s}{n}]\neq\emptyset$ we choose
$s_n\in[\frac{s}{n}]$, and we write
$\Lambda(S):=\{s_n:s\in S, 0\neq n\in\ZZ\}\cup\{0\}$.
Then $|\Lambda(S)|\leq|S|\cdot\aleph_0$, and $\mbox{div}(S)=\Lambda(S)+t(K)$.

\begin{lemma} 
Let $M$ be an Abelian group such that $|M|=\kappa>\omega$ and
let $S\in[M]^{<\kappa}$, $E\in[M]^\kappa$ with $S$ a subgroup. If either
\begin{itemize}
\item[(i)] $t(M)<\kappa$, or 
\item[(ii)] there is $p\in\PP$ such that $p\cdot M=\{0\}$,
\end{itemize}
 then there is $x\in E$ such that $\langle x\rangle\cap S=\{0\}$. In
case {\rm (i)}, $x$ may be chosen in $M\backslash t(M)$.
\end{lemma}

\begin{proof} (i) From $\mbox{div}(S)=\Lambda(S)+t(M)$ follows $|\mbox{div}(S)|<\kappa$, and any
$x\in E\backslash\mbox{div}(S)\subseteq E\backslash t(M)$ is as required.

(ii) Since
$M\simeq\bigoplus_\kappa\,\ZZ(p)=\bigoplus_{\xi<\kappa}\,\ZZ(p)_\xi$, 
there 
is $A\in[\kappa]^{<\kappa}$ such that $S\subseteq\bigoplus_{\xi\in
A}\,\ZZ(p)_\xi$. Any $x\in E$ such that $0\neq x\notin\bigoplus_{\xi\in
A}\,\ZZ(p)_\xi$ is as required.
\end{proof}

\begin{theorem}  
Let $M$ be an infinite, compact, metrizable, Abelian group such that
either
\begin{itemize}
\item[(i)] $|t(M)|<\cc$ or
\item[(ii)] there is $p\in\PP$ such that $p\cdot M=\{0\}$.
\end{itemize}
Then $M$ admits a dense, nonmeasurable subgroup $D$ such that
$|M/D|=\cc$.

In case {\rm (i)} one may arrange $D\simeq\bigoplus_\cc\,\ZZ$, in case
{\rm (ii)} one may arrange $D\simeq\bigoplus_\cc\,\ZZ(p)$.
\end{theorem}

\begin{proof}
Let $\{F_\xi:\xi<\cc\}$ be an enumeration of all uncountable,
closed subsets of $M$, and define $E_\xi:=(F_\xi\backslash
t(M))\backslash\{0\}$ in case (i), $E_\xi:=F_\xi\backslash\{0\}$ in case
(ii). It is a theorem of Cantor~\cite{cantor84}(page~488)
that each $|F_\xi|=\cc$
(see \cite{haus14}(VIII\,\S9\,II) or \cite{engel}(4.5.5(b)) for more
modern treatments); hence each $|E_\xi|=\cc$.
There is $x_0\in E_0$,
and by Lemma~2.3 there is $y_0\in E_0$ such that $\langle x_0\rangle\cap\langle
y_0\rangle=\{0\}$. 

Now let $\xi<\cc$, suppose that $x_\eta$, $y_\eta$ have been chosen
for all $\eta<\xi$, and apply Lemma~2.3 twice to choose $x_\xi,y_\xi\in
E_\xi$ such that
\begin{align*}
\langle\{x_\xi\}\rangle\cap\langle\{x_\eta:\eta<\xi\}\cup\{y_\eta:\eta<\xi\}\rangle&
=\{0\},\ \text{and}\\
\langle\{y_\xi\}\rangle\cap\langle\{x_\eta:\eta\leq\xi\}\cup\{y_\eta:\eta<\xi\}\rangle&
=\{0\}.
\end{align*}
Thus $x_\xi,y_\xi$ are defined for all $\xi<\cc$. We define
$D:=\langle\{x_\xi:\xi<\cc\}\rangle$. Clearly $D=\bigoplus_{\xi<\cc}\,\ZZ$
in case (i), and $D=\bigoplus_{\xi<\cc}\,\ZZ(p)_\xi$ in case (ii) since
$|D|=\cc$ and $p\cdot D=\{0\}$. For $\xi<\eta<\cc$ we have
$y_\eta+D\neq y_\xi+D$, so $\cc\geq|K/D|\geq\cc$.

For nonempty open $U\subseteq K$  there is by the regularity of
$\lambda$ a (necessarily uncountable) compact set $F=F_\xi\subseteq U$
such that $\lambda(F_\xi)>0$. Then $x_\xi\in E_\xi\cap D\subseteq F_\xi\cap
D$. Thus $D$ is dense in $K$. If $D\in\sM(K)$ with $\lambda(D)>0$
then $D=K$ by Theorem~1.4(c), contrary
to the relation $|K/D|=\cc$. If $D\in\sM(K)$ with $\lambda(D)=0$
then $\lambda(K\backslash D)=1$ and there is $F_\xi\subseteq K\backslash
D$ such that $\lambda(F_\xi)>0$; then $y_\xi\in E_\xi\cap D\subseteq
D\backslash D=\emptyset$, a contradiction.
\end{proof}

\begin{corollary}
Let $M$ be a (compact, Abelian, metrizable) group of one of
these types.
\begin{itemize}
\item[(i)] $M=\TT$;
\item[(ii)] $M=\Delta_p$ ($p\in\PP$), the group of $p$-adic integers;
\item[(iii)] $M=\prod_{k<\omega}\,\ZZ(p_k)$, $p_k\in\PP$, $(p_k)_k$ 
faithfully indexed;
\item[(iv)] $M=(\ZZ(p))^\omega$ ($p\in\PP$).
\end{itemize}
Then $M$ admits a dense, nonmeasurable subgroup $D$ such that $|M/D|=\cc$.
\end{corollary}

\begin{proof} Surely $|t(\TT)|=\omega$, and  
$t(\prod_{k<\omega}\,\ZZ(p_k))$ is the countable group\break 
$\bigoplus_{k<\omega}\,\ZZ(p_k)$. 
If $0\neq h\in\Delta_p=\Hom(\ZZ(p^\infty),\TT)$, then
$h[\ZZ(p^\infty)]\simeq\break\ZZ(p^\infty)/\ker(h)\simeq\ZZ(p^\infty)$ 
since
$|\ker(h)|<\omega$, so $h[\ZZ(p^\infty)]$ is not of bounded order; thus
$t(\Delta_p)=\{0\}$. It follows for $M$ as in (i), (ii) and (iii) that
$|t(M)|=1<\cc$ or $|t(M)|=\omega<\cc$, so
Theorem~2.4(i) applies.
For $M$ as in (iv) surely $p\cdot M=\{0\}$, so Theorem~2.4(ii)
applies.
\end{proof} 

\begin{theorem}   
Let $K$ be an infinite, compact, Abelian metrizable group. Then
($|K|=\cc$ and) $K$ has a dense, nonmeasurable subgroup $H$ such that $|H|=\cc$ and
$|K/H|=\cc$.
\end{theorem}

\begin{proof}
 The (discrete) dual group $G=\widehat{K}$ satisfies
$|G|=w(K)=\omega$. As with any countably infinite Abelian group, $G$
must satisfy (at least) one of these conditions: (i)~$r_0(G)>0$:
(ii)~$|G_p|=\omega$ with $r_p(G)<\omega$ for some $p\in\PP$;
(iii)~$0<r_p(G)<\omega$ for infinitely many $p\in\PP$;
(iv)~$r_p(G)=\omega$ for some $p\in\PP$. According as (i), (ii), (iii)
or (iv) holds we have, respectively, $G\supseteq\ZZ$,
$G\supseteq\ZZ(p^\infty)$, $G\supseteq\bigoplus_{k<\omega}\,\ZZ(p_k)$, or
$G\supseteq\bigoplus_\omega\,\ZZ(p)$, so taking adjoints we have a
continuous surjection $\phi$ from $K$ onto a group $M$ of the form
$\widehat{\ZZ} = \TT$,
$\widehat{\ZZ(p^\infty)} = \Delta_p$,
$\widehat{\bigoplus_{k<\omega}\,\ZZ(p_k)}
= \prod_{k<\omega}\ZZ(p_k)$, or
$\widehat{\bigoplus_\omega\,\ZZ(p)} = (\ZZ(p))^\omega$. 
According to Corollary~2.5 the group $M$ has a dense, nonmeasurable subgroup $D$
such that $|M/D|=\cc$, and then by Theorem~2.2 with $H:=\phi^{-1}(D)$ the group
$H$ is dense and nonmeasurable in $K$. If $a,b\in K$ with $a+H=b+H$ then
$\phi(a)+D=\phi(b)+D$, so $\cc=|K|\geq|K/H|\geq|M/D|=\cc$.
\end{proof}

\begin{theorem} 
Let $K$ be an infinite, compact, Abelian, metrizable group.
Then $K$ admits a family of $2^{|K|}$-many dense, nonmeasurable subgroups,
each of cardinality $\cc$.
\end{theorem}

\begin{proof}
Let $H$ be as given in Theorem~2.6 and let
$$\sH:=\{S:H\subseteq S\subseteq K,\ \text{$S$ is a proper subgroup of 
$K$}\}.$$
Then $|\sH|=2^\cc$ by Theorem~1.7. Theorem~1.4(c) shows for
$S\in\sH$ that $S\in\sM(K)$ with $\lambda(S)>0$ is impossible, and if
$S\in\sH$ with $\lambda(S)=0$ then $\lambda(H)=0$, a
contradiction.
\end{proof}

\begin{remark}
{\rm
For our application in Theorem~4.1 below we do not require
that $|K/S|=\cc$, but in fact that condition does hold for
$2^\cc$-many $S\in\sH$.
}
\end{remark}

\section{Many Nonmeasurable Subgroups: The Nonmetrizable Case}

We turn now to the case $w(K)=\kappa>\omega$. Again, our
goal is to show that such a compact Abelian group $K$ contains
$2^{|K|}$-many dense, nonmeasurable subgroups of cardinality $|K|$.
By a well-known structure
theorem (see Theorem~3.2 below), such a group $K$ admits a continuous, surjective
homomorphism onto a group $M$ of the form $M=\prod_{\xi<\kappa}\,M_\xi$ 
with
each $M_\xi$ a (compact) subgroup of $\TT$, and according to Theorem~2.2 it suffices to
show that such $M$ has a family of $2^{|M|}$-many ($=2^{|K|}$-many) such
subgroups. We find it convenient to show a bit more, namely that $M$, and hence also $K$, admits a
family of $2^{|M|}$-many subgroups each of which is $G_\delta$-dense in $M$.
In the transition, we will invoke the following lemma.

\begin{lemma}  
A proper, $G_\delta$-dense subgroup $H$ of a compact group $K$ is
nonmeasurable.
\end{lemma}

\begin{proof}
As usual, using Theorem~1.4(c), $H\in\sM(K)$ with
$\lambda(H)>0$ is impossible; it suffices then to show that
$\lambda(H)=0$ is also impossible. The following argument is from
\cite{kodaira41} and \cite{kakkod}, as exposed by Halmos~\cite{Ha}.
If $\lambda(K\backslash H)=1>0$ there are a compact set $C$ and a Baire
set $F$ of $K$ such that $F\subseteq C\subseteq K\backslash H$ and
$\lambda(F)=\lambda(C)>0$ (\cite{Ha}(64H and p.~230)). As with any
nonempty Baire set, $F$ has the form $F=XB$ for a suitably chosen
compact Baire subgroup $B$ of $K$ and $X\subseteq K$ (\cite{Ha}(64E)).
Since every compact Baire set is a $G_\delta$-set (\cite{Ha}(51D)), each
$x\in X$ has $xB$ a $G_\delta$-set. From $xB\cap H=\emptyset$ it follows
that $H$ is not $G_\delta$-dense in $K$, a contradiction.
\end{proof}

Now for compact Abelian $K$ with $w(K)=\kappa>\omega$ we write
$G=\widehat{K}$, so that $|G|=w(K)=\kappa$, and we denote the
torsion-free rank and (for $p\in\PP$) the $p$-rank of $G$ by $\kappa_0=r_0(G)$
and $\kappa_p=r_p(G)$, respectively. Since $w(K)=|G|>\omega$ we have
$|G|=\kappa=\kappa_0+\sum_{p\in\PP}\,\kappa_p$
(with perhaps $\kappa_i=0$ for certain $i\in\PP\cup\{0\}$), and algebraically
\begin{equation*}
G\supseteq\bigoplus_{\kappa_0}\,\ZZ \oplus \bigoplus_{p\in\PP}
\bigoplus_{\kappa_p}\,\ZZ(p).
\tag{$\ast$}
\end{equation*}

\begin{theorem}  
Let $K$ be a compact, Abelian group such that $w(K)=\kappa>\omega$.
Then $K$ has a family of $2^{|K|}$-many $G_\delta$-dense subgroups of
cardinality $|K|$.
\end{theorem}

\begin{proof}
 The map $\psi$ adjoint to the inclusion map in ($\ast$) is a
continuous, surjective homomorphism:
$$\psi:K
= 
\widehat{G}\twoheadrightarrow\TT^{\kappa_0}\times\prod_{p\in\PP}\,(\ZZ(p))^{\kappa_p}$$
(with, again, perhaps $\kappa_i=0$ for certain
$i\in\PP\cup\{0\}$). If $\kappa_0=\kappa$ or some $\kappa_p=\kappa$
($p\in\PP$)---as is necessarily the case in $\cf(\kappa)>\omega$---then a
suitable projection map $\pi$ furnishes a continuous, surjective
homomorphism $\phi=\pi\circ\psi:K\twoheadrightarrow M:=F^\kappa$ with
$F=\TT$ or $F=\ZZ(p)$. That such a compact group $M$ admits a faithfully
indexed family $\{D_\xi:\xi<2^\kappa\}$ of $G_\delta$-dense subgroups is known (see
\cite{wclr:88}(4.4)); according to \cite{wclr:82}(2.2) the groups
$\phi^{-1}(D_\xi)$ ($\xi<2^\kappa$) are then $G_\delta$-dense in $K$.
We assume now that $\kappa_0<\kappa$ and each $\kappa_p<\kappa$
($p\in\PP$), so that necessarily $\cf(\kappa)=\omega$; we show that the
``ultrafilter technique" of \cite{wclr:82} and \cite{wclr:88} can
be adapted to cover this case also.

The set of uniform ultrafilters over an (infinite, discrete) set $A$ is
here denoted $\uu(A)$. It is well known that $|\uu(A)|=2^{2^{|A|}}$.
(See~\cite{WCSN:74} for a proof and for other relevant combinatorial facts.)

Let $(p_n)_n$ be a faithfully indexed sequence in $\PP$ such
that $\omega<\kappa_{p_n}<\kappa$ and
$\sup_{n<\omega}\,\kappa_{p_n}=\kappa$, and choose pairwise disjoint sets
$A_n$ such that $|A_n|=\kappa_{p_n}$; then with
$A:=\bigcup_{n<\omega}\,A_n$ we have $|A|=\kappa$. We write
$\phi=\pi\circ\psi:K\twoheadrightarrow
M=\prod_{n<\omega}\,\ZZ(p_n)^{A_n}$, we view $A$ as a discrete space
and we view each $f\in M$ as a function of the form
$f=\bigcup_{n<\omega}\,f_n$ with $f_n:A_n\rightarrow\ZZ(p_n)\subseteq\TT$.
The Stone-\v{C}ech extension of $f$, mapping $\beta(A)$ into $\TT$,
is denoted $\overline{f}$. Now for each $q\in\uu(A)\subseteq\beta(A)$
we define $h_q:M\rightarrow\TT$ by
$h_q(f)=\overline{f}(q)\in\TT$. For $f_0,f_1\in M$ we have
$\overline{f_0+f_1}|A=f_0+f_1=\overline{f_0}|A+\overline{f_1}|A$, so
$\overline{f_0+f_1}=\overline{f_0}+\overline{f_1}$; thus
$h_q\in\Hom(M,\TT)$ and $\graph(h_q)$ is a subgroup of $M\times\TT$.

We claim that $\graph(h_q)$ is $G_\delta$-dense in $M\times\TT$. (It
will follow in particular then that $M_q:=\ker(h_q)$ is a
$G_\delta$-dense subgroup of $M$.)

Let $C=\{a_k\}_k$ be a faithfully indexed subset of $A$, say with $a_k\in
A_{n_k}$, and let $v_k\in\ZZ(p_{n_k})$; further, let $w\in\TT$. Then
with $V:=\{f\in M:f(a_k)=v_k\}$ and $W:=\{w\}$, the set $V\times W$
is a nonempty $G_\delta$-subset of $M\times\TT$; since
each nonempty $G_\delta$-subset of
$M\times\TT$ contains a set of the form $V\times W$, to prove the claim
it suffices to show that for each $q\in\uu(A)$ there is
$f\in M$ such that $f(a_k)=v_k$ (that is, $f\in V$) and also
$\overline{f}(q)=w$ (that is,  $h_q(f)\in W$). Since
$\bigcup_{n<\omega}\,\ZZ(p_n)$ is dense in $\TT$ there is
$w_n\in\ZZ(p_n)$ such that $w_n\rightarrow w$. We define $g\in M$ by $g_n\equiv
w_n$ on $A_n$, and we define $f\in M$ by
$$f(a)
= \left\{
\begin{array}{cl} 
g(a) & \mbox{if $a\in A\backslash C$,}\\
w_{n_k}   & \mbox{if $a = a_k\in C$.}
\end{array}\right.
$$
Each neighborhood $U$ of $q$ in $\beta(A)$ meets infinitely many of
the sets $A_n$, so $w_n\in g[U]$ for infinitely many $n<\omega$ and
hence $\overline{g}(q)=w$. Clearly $f\in V$, and since $f\equiv g$ on
$A\backslash C\in q$ we have
$h_q(f)=\overline{f}(q)=\overline{g}(q)=w$. The claim is proved.

Since $|\uu(A)|=2^{2^\kappa}$, to see that $M$
has $2^{|M|}$-many $G_\delta$-dense subgroups it
suffices to show that if $q_0\neq q_1$ with $q_i\in\uu(A)$ then
$M_{q_0}=\ker(h_{q_0})\neq\ker(h_{q_1})=M_{q_1}$. Choose $F\subseteq A$
such that $F\in q_0$ and $A\backslash F\in q_1$, fix (for some $n<\omega$)
$t\in\ZZ(p_n)\cap(\frac{1}{4},\frac{3}{4})\subseteq\TT$ and
define $f\in M$ by
$$f(a)
= \left\{
\begin{array}{cl} 
0 &\mbox{if $a\in F$},\\
t & \mbox{if $a\in A \setminus F$.}
\end{array}\right.
$$
Then $h_{q_0}(f)=\overline{f}(q_0)=0$ and
$h_{q_1}(f)=\overline{f}(q_1)=t\neq0$, so $f\in\ker(h_{q_0})=M_{q_0}$ while
$f\notin\ker(h_{q_1})=M_{q_1}$.

Again as in \cite{wclr:82}(2.2), each $\phi^{-1}(M_q)$
($q\in\uu(A)$) is a (proper) $G_\delta$-dense subgroup of $K$.

It remains to show that the $G_\delta$-dense subgroups of $K$
constructed above are of cardinality $|K|$. If $|K|=2^\kappa=\cc$ (as
can occur for certain $\kappa>\omega$ in some models of set theory),
this is clear since a $G_\delta$-dense subgroup of a compact group is
pseudocompact \cite{wckr:66} and hence, if infinite, has cardinality at
least $\cc$~\cite{vdie},~\cite{wclr:85}. If $|K|=|M|=2^\kappa>\cc$ then from
$h_q[M]\simeq M/M_q$ and $|M/M_q|\leq|\TT|=\cc$ we infer $|M_q|=|M|$ and
hence $|K|\geq|\phi^{-1}(M_q)|\geq|M_q|=|M|=|K|$. The same argument applies to
the $G_\delta$-dense subgroups $D_\xi$ of $F^\kappa$ given in
\cite{wclr:88}(4.4) which appear in the first paragraph of this proof,
for each of those is the kernel of a homomorphism from $F^\kappa$ onto
$F\subseteq\TT$.
\end{proof}

We note in passing that in general not every $G_\delta$-dense subgroup
$D$ of a compact, nonmetrizable group $K$ satisfies $|D|=|K|$. See in
this connection Remark~6.4 below.

Corollary 3.3 is now immediate from Lemma~3.1 and Theorem~3.2; and
Theorem~3.4, which is (1) of our Abstract, is the conjuction of
Theorem~2.7 and Corollary~3.3.

\begin{corollary}  
Let $K$ be a compact, Abelian group such that $w(K)>\omega$.
Then $K$ admits a family of $2^{|K|}$-many dense, nonmeasurable subgroups,
each of cardinality $|K|$.
\end{corollary}

\begin{theorem} 
Every infinite, Abelian compact (Hausdorff) group $K$ admits
$2^{|K|}$-many dense, non-Haar-measurable subgroups of cardinality
$|K|$. When $K$ is nonmetrizable, these may be chosen to be
pseudocompact. 
\end{theorem}

\begin{remark} 
{\em It is known \cite{Str92} that $\RR$ contains
a dense non-measurable subgroup of both cardinality and index $\cc$. 
Arguing as in Theorem 2.7, we see then that $\RR$ has
$2^\cc$-many dense non-measurable subgroups of both cardinality and index 
$\cc$. We say as usual that a {\it compactly generated}
group is a topological group
generated, in the algebraic sense, by a compact subset. By 
\cite{ehkr:63} (9.8) for every (Hausdorff) locally compact Abelian 
compactly generated group $G$ there are non-negative integers $m$ 
and $n$ and a compact Abelian group $K$ such that $G$ is of the 
form $\RR^m \times K \times \ZZ^n$. If $G$ is not discrete, then 
either $m >0$ or $K$ is not discrete, so $w(G) = \omega + w(K)$. It
follows  that a non-discrete (Hausdorff) locally 
compact Abelian compactly generated group $G$ has $2^{|G|}$-many 
dense non-Haar-measurable subgroups of both cardinality and index $|G|$. 
More generally, let $G$ be an arbitrary non-discrete (Hausdorff) locally 
compact Abelian group, and let $H$ be an open compactly generated 
subgroup of $G$. It follows that $G$ has $2^{|H|}$-many dense 
non-Haar-measurable subgroups $\{G_\xi: \xi < 2^{|H|} \}$,
each of cardinality $|G|$, and such that each $H_\xi:=G_\xi \cap H$
has $|H/H_\xi|=|H|$. We leave the details to the reader.}
\end{remark}

\section{Group Topologies Without Convergent Sequences}

Here we pull together the threads of Sections 2 and 3.

\begin{theorem}  
Every infinite Abelian group $G$ admits a family $\sA$ of
totally bounded group topologies, with $|\sA|=2^{2^{|G|}}$, such that no 
nontrivial sequence in $G$ converges in any of the topologies in $\sA$. 
One may arrange in addition that
\begin{itemize}
\item[(i)] $w(G,\sT)=2^{|G|}$ for each $\sT\in\sA$, and
\item[(ii)] for distinct $\sT_0$, $\sT_1\in\sA$ the spaces $(G,\sT_0)$ 
and $(G,\sT_1)$ are not homeomorphic.
\end{itemize}
\end{theorem}

\begin{proof}
By Theorem 2.7 when $|G|=\omega$, and by Corollary~3.3 when
$|G|>\omega$, the compact group $K:=\Hom(G,\TT)=\widehat{G_d}$ admits a
family $\sH$ of dense, nonmeasurable subgroups such that $|\sH|=2^{2^{|G|}}$ and
$|H|=|K|$ for each $H\in\sH$. According to Theorems~1.5 and 1.6 the family
$\sA:=\{\sT_H:H\in\sH\}$  satisfies all requirements except (perhaps)
(ii). A homeomorphism between two of the spaces $(G,\sT_0)$, $(G,\sT_1)$
with $\sT_i\in\sA$ is realized by a permutation of $G$, and there are
just $2^{|G|}$-many such functions, so for each $\sT\in\sA$ there are at
most $2^{|G|}$-many $\sT'\in\sA$ such that $(G,\sT)=_h(G,\sT')$.
Statement (ii) then follows (with
$\sA$ replaced if necessary by a suitably chosen subfamily of
cardinality $2^{|K|}=2^{2^{|G|}}$).
\end{proof}

\begin{remark} 
{\rm
The case $G=\ZZ$ of Theorem~4.1 is not new. See in this connection
\cite{raczphd} and \cite{racz:02}, which were the motivation for much of
the present paper.
}
\end{remark}

\section{Topologies With Convergent Sequences}

We turn now to the complementary or opposing problem, that of finding on
an arbitrary infinite Abelian group $G$ the maximal number (that is,
$2^{2^{|G|}}$-many) totally bounded group topologies in which some nontrivial sequence
converges. Again, this result was first achieved for $G=\ZZ$ in
\cite{raczphd}, \cite{racz:02}: there, $2^\cc$ many topologies in $\ttt(\ZZ)$
are constructed in which
an arbitrary sequence $x_n$, fixed in advance and
satisfying $x_{n+1}/x_n\geq n+1$, converges to $0$. The condition
$x_{n+1}/x_n\geq n+1$ of \cite{raczphd} and \cite{racz:02}
was relaxed in \cite{Barb} to
$x_{n+1}/x_n\to\infty$.

To handle the case of general Abelian $G$, our
strategy is to show first that certain ``basic" countable groups
accept $2^\cc$ many such
topologies. We begin with technical results concerning groups of the
form $\bigoplus_{k<\omega}\,\ZZ(p_k^{r_k})$ and of the form 
$\ZZ(p^\infty)$.

We remark for emphasis that in Theorem~5.1 the given sequence $(p_k)_k$ in
$\PP$ is not necessarily faithfully indexed. Indeed the case
$p_k=p\in\PP$ (a constant sequence) is not excluded. For $x\in
A=\bigoplus_{k<\omega}\,\ZZ(p_k^{r_k})$ we write $x=(x(k))_{k<\omega}$.

\begin{theorem}   
Let $p_k\in\PP$ and $A=\bigoplus_{k<\omega}\,\ZZ(p_k^{r_k})$
with $0<r_k<\omega$, and let $(x_n)_{n<\omega}$ be a faithfully indexed
sequence in $A$ such that
\begin{itemize}
\item[(i)] there is $S\in[\omega]^\omega$ such that $x_n(k)=0$ for all 
$n<\omega$, $k\in S$, and
\item[(ii)] $|\{n<\omega:x_n(k)\neq0\}|<\omega$ for all $k<\omega$.
\end{itemize}
Let $\{A_\xi:\xi<\cc\}$ enumerate
$\sP(S)\cup[\omega]^{<\omega}$, and for $\xi<\cc$ define
$h_\xi\in\Hom(A,\TT)$ by
$h_\xi(x)=\sum_{k\in A_\xi}\,x(k)$.
Then
\begin{itemize}
\item[(a)] the set $\{h_\xi:\xi<\cc\}$ is faithfully indexed;
\item[(b)] the set $\{h_\xi:\xi<\cc\}$ separates points of $A$; and
\item[(c)] $h_\xi(x_n)\rightarrow0$ for each $\xi<\cc$.
\end{itemize}
\end{theorem}

\begin{proof}
If $\xi,\xi'<\cc$ with $\xi\neq\xi'$, say
$k\in A_\xi\backslash A_{\xi'}$, then any $x\in G$
such that $0\neq x(k)\in\ZZ(p_k^{r_k})$ and $x(m)=0$ for $k\neq m<\omega$
satisfies
$$h_\xi(x)=x(k)\neq0=h_{\xi'}(x).$$

(b) Let $x,x'\in G$ with, say, $x(k)\neq x'(k)$, and let $\{k\}=A_\xi$
($\xi<\cc$). Then
$$h_\xi(x)=x(k)\neq x'(k)=h_\xi(x').$$

(c) If $A_\xi\in\sP(S)$ then $h_\xi(x_n)=0$ for all $n$ by (i), so
$h_\xi(x_n)\rightarrow0$. If $A_\xi\in[\omega]^{<\omega}$ then by (ii)
there is $N<\omega$ such that $h_\xi(x_n)=0$ for all $n>N$, so again
$h_\xi(x_n)\rightarrow0$.
\end{proof} 

Next, following \cite{ehkr:63} and \cite{fuch:70}, we identify the 
elements of the compact group $\Delta_p=\Hom(\ZZ(p^\infty),\TT)$ with 
those sequences $h=(h(k))_{k<\omega}$ of integers such that $0\leq 
h(k)<p-1$ for all $k<\omega$. For $\frac{a}{p^n}\in\ZZ(p^\infty)$ (with 
$0\leq a<p^n-1$)
we have
$$h(\frac{a}{p^n})
= a \cdot h(\frac{1}{p^n})
= a \cdot \sum_{k=0}^{n-1}\,\frac{h(k)}{p^{n-k}}
= \frac{a}{p^n} \cdot \sum_{k=0}^{n-1}\,h(k) \cdot p^k\quad 
(\mod~1).$$
In what follows we write $\Fac:=\{n!:n<\omega\}$.

\begin{theorem}  
Let $(a_n)_{n<\omega}$ be a sequence of integers such that
$0\leq a_n<p-1$ for all $n<\omega$, and let
$x_n=\frac{a_n}{p^{n!}}\in\ZZ(p^\infty)$. Let $\{A_\xi:\xi<\cc\}$
enumerate $\sP(\Fac)\cup[\omega]^{<\omega}$, and for $\xi<\cc$ define
$h_\xi=(h_\xi(k))_{k<\omega}\in\Delta_p$ by
$$h_\xi(k)=\left\{
\begin{array}{lc}
1 &\mbox{if $k\in A_\xi$,}\\
0 &\mbox{otherwise}.
\end{array}\right.
$$
Then
\begin{itemize}
\item[(a)] the set $\{h_\xi:\xi<\cc\}$ is faithfully indexed;
\item[(b)] the set $\{h_\xi:\xi<\cc\}$ separates points of $A$; and
\item[(c)] $h_\xi(x_n)\rightarrow0$ for each $\xi<\cc$.
\end{itemize}
\end{theorem}

\begin{proof}
The proofs of (a) and (b) closely parallel their analogues in
Theorem~5.1. We prove (c). If $A_\xi\in[\omega]^{<\omega}$ there is
$N<\omega$ such that
$$h_\xi(x_n)
= \frac{a_n}{p^{n!}} 
\cdot \sum_{k=0}^{n!-1}\,h_\xi(k)
\cdot p^k\leq\frac{p-1}{p^{n!}}
\cdot \sum_{k=0}^N\,p^k\quad\text{for all $n>N$,}$$
so $h_\xi(x_n)\rightarrow0$. If $A_\xi\in\sP(\Fac)$ then
\begin{align*} 
h_\xi(x_n)&
= \frac{a_n}{p^{n!}}
\cdot \sum_{k=0}^{n!-1}\,h_\xi(k)
\cdot p^k\\
&
= \frac{a_n}{p^{n!}}
\cdot \sum_{k\in\Fac,0 \leq k<n!-1} h_\xi(k)
\cdot p^k\\
&
\leq \frac{p-1}{p^{n!}}
\cdot \sum_{m=0}^{n-1}\,p^{m!}\\ 
& 
\leq \frac{p-1}{p^{n!}}
\cdot n
\cdot p^{(n-1)!}\\
&
< \frac{(p-1) \cdot n}{p^n},
\end{align*}
and again $h_\xi(x_n)\rightarrow0$.
\end{proof}

\begin{theorem}  
Let $A=\ZZ$ or $A=\bigoplus_{k<\omega}\,\ZZ(p_k)$ ($p_k\in\PP$,
repetitions allowed), or $A=\ZZ(p^\infty)$ ($p\in\PP$). Then there are a
faithfully indexed sequence $(x_n)_{n<\omega}$ in $A$ and a point-separating subgroup
$H\subseteq\Hom(A,\TT)$ such that $|H|=\cc$ and $x_n\rightarrow0$ in the
space $(A,\sT_H)$.
\end{theorem}

\begin{proof}
We refer the reader for a detailed proof in the case $A=\ZZ$ to
\cite{raczphd},~\cite{racz:02}. When $A=\bigoplus_{k<\omega}\,\ZZ(p_k)$ or
$A=\ZZ(p^\infty)$ there are, according to Theorem~5.1 or 5.2
respectively, a  sequence $(x_n)$ in $A$ and a faithfully indexed family
$\{h_\xi:\xi<\cc\}\subseteq\Hom(A,\TT)$ such that $h_\xi(x_n)\rightarrow0$ for each
$\xi<\cc$. It is then clear, as in \cite{CTWu:93}, that with
$H:=\langle\{h_\xi:\xi<\cc\}\rangle\subseteq\Hom(A,\TT)$ we have
$h(x_n)\rightarrow0$ for each $h\in H$.
\end{proof}

\begin{corollary}
Let $G$ be an infinite Abelian group. There are a faithfully indexed
sequence $(x_n)_{n<\omega}$ in $G$ and a topology $\sT_{H^*}\in\ttt(G)$
such that $|H^*|=2^{|G|}$ and $x_n\rightarrow0$ in $(G,\sT_{H^*})$.
\end{corollary}

\begin{proof}
As indicated earlier, $G$ contains algebraically a group $A$ such that
$A\simeq\ZZ$ or $A\simeq\bigoplus_{k<\omega}\,\ZZ(p_k)$ or 
$A\simeq\ZZ(p^\infty)$.
Let $H\subseteq\Hom(A,\TT)$ and $(x_n)_n$ in $A$
be as in Theorem~5.3, and set
$H^*:=\{k\in\Hom(G,\TT):k|A\in H\}$. Clearly $k(x_n)\rightarrow0$ for
each $k\in H^*$, so $x_n\rightarrow0$ in $(G,\sT_{H^*})$. If $|G|=\omega$
then $|H^*|=\cc=2^\omega$ since $\cc=|H|\leq|H^*|\leq|\TT^G|=\cc$.  If
$|G|>\omega$ we write
$$\AA(\widehat{G},A)
:= \{k\in\Hom(G,\TT):k\equiv0\ \text{on $A$}\}\subseteq H^*;$$
then $|\AA(\widehat{G},A)|=2^{|G/A|}=2^{|G|}$ since algebraically
$\AA(\widehat{G},A)=\Hom(G/A,\TT)$ and $|G|=|G/A|$, so $|H^*|=2^{|G|}$ in this
case also.
\end{proof}

We have arrived at the final result of this Section, which we view as a
companion or ``echo" to Theorem~4.1.

\begin{theorem}
Every infinite Abelian group $G$ admits a family $\sB$ of totally bounded 
group topologies, with $|\sB|=2^{2^{|G|}}$, such that some (fixed) nontrivial
sequence in $G$ converges in each of the topologies in $\sB$. One may 
arrange in addition that
\begin{itemize}
\item[(i)]
$w(G,\sT)=2^{|G|}$ for each $\sT\in\sB$, and
\item[(ii)] for distinct $\sT_0$, $\sT_1\in\sB$ the spaces $(G,\sT_0)$ 
and $(G,\sT_1)$ are not homeomorphic.
\end{itemize}
\end{theorem}

\begin{proof}
Let $H^*$ be a subgroup of $\Hom(G,\TT)$ such that $|H^*|=2^{|G|}$ and some
nontrivial sequence $(x_n)_n$ in $G$ satisfies $x_n\rightarrow0$ in
$(G,\sT_{H^*})$. There is a subgroup $H$ of $H^*$ such
that $H$ separates points of $G$
and $|H|=|G|$, and since $|H^*/H|=2^{|G|}$ there is by Theorem~1.7
a faithfully indexed family
$\{H_\xi:\xi<2^{2^{|G|}}\}$ of (point-separating) groups such that $H\subseteq
H_\xi\subseteq H^*$ for each $\xi$. Then
$\sB:=\{\sT_{H_\xi}:\xi<2^{2^{|G|}}\}$ satisfies all requirements except
perhaps (ii), and (ii) is handled as in the final sentences of the proof
of Theorem~4.1.
\end{proof}

\begin{remarks}
{\rm
(a) If $\sA$ and $\sB$ are as in Theorems~4.1 and 5.5, and if
$\sT_0,\sT_1\in\sA\cup\sB$ with $\sT_0\neq\sT_1$, then the spaces
$(G,\sT_0)$ and $(G,\sT_1)$ are not homeomorphic. For if both $\sT_i\in\sA$
or both $\sT_i\in\sB$ this is already proved, while if (say)
$\sT_0\in\sA$ and $\sT_1\in\sB$ then $(G,\sT_1)$ has a nontrivial
convergent sequence and $(G,\sT_0)$ does not.

(b) We emphasize that for an infinite Abelian group $G$, the algebraic
structure of a point-separating subgroup $H\subseteq\Hom(G,\TT)$ by no
means determines the topology $\sT_H$ on $G$. It is noted explicitly in
\cite{raczphd}, \cite{racz:02} that when $G=\ZZ$ then every one of the
topologies in the families $\sA$ and $\sB$ (as in Theorems~4.1 and 5.5)
can be chosen of the form $\sT_H$ with $H\subseteq\TT=\Hom(G,\TT)$ and
with $H\simeq\bigoplus_{\xi<\cc}\,\ZZ_\xi$.
}
\end{remarks}

\section{Concluding Remarks}

It was shown in \cite{Barb}, assuming MA, that there is a measurable subgroup
$H$ of $\TT$ such that $\lambda(H)=0$ and no nontrivial sequence converges in
the space $(\ZZ,\sT_H)$. Responding to a question in \cite{Barb}, the
authors of \cite{HaKu} achieved the same result in ZFC. It is natural to
ask if the comparable statement holds for each infinite Abelian group.
In detail:

\begin{question}  
{\rm
Does every infinite Abelian group $G$ admit a group topology
$\sT=\sT_H\in\ttt(G)$ with no nontrivial convergent sequences such that $H$ is
measurable in the compact group $\Hom(G,\TT)=\widehat{G_d}$ and
$\lambda(H)=0$?
}
\end{question}

Our convergent sequences in Theorems~5.1 and 5.2 were quite ``thin".
Viewing those results from another perspective, a natural
question arises.

\begin{question}   
{\rm
Let $G$ be an infinite Abelian group. For which faithfully indexed
sequences $(x_n)_n$ in $G$ is there a topology $\sT_H\in\ttt(G)$ with
$w(G,\sT_H)=|H|=2^{|G|}$ such that $x_n\rightarrow0$ in $(G,\sT_H)$?
}
\end{question}

Many questions arise in the non-Abelian context. Perhaps this one of
Saeki and Stromberg~\cite{SaStr:85} bears repeating.

\begin{question}
\cite{SaStr:85}. {\rm Does every infinite (not necessarily
Abelian) compact group have a dense, nonmeasurable subgroup?}
\end{question}

\begin{remark}
{\rm
Malykhin and Shapiro~\cite{malyshap} showed by a direct argument that
for every faithfully indexed sequence $x=(x_n)_n$ in an Abelian group $G$ there is
$h_x\in\Hom(G,\TT)$ such that $h_x(x_n)\not\rightarrow0$. Thus every
Abelian $G$ with $|G|=\alpha$ admits a topology $\sT_H\in\ttt(G)$ such that
$|H|=w(G,\sT_H)=\alpha^\omega$ and no nontrivial sequence converges in
$(G,\sT_H)$. This statement can be improved slightly using Theorem~1.6
above and this result from \cite{wclr:85}(2.2): Every compact group $K$
with $w(K)=\alpha$ contains a dense, countably compact, subgroup $H$ such that
$|H|=(\log\,\alpha)^\omega$. Since such $H$ is pseudocompact and hence
nonmeasurable, we have (beginning with Abelian $G$ such that $|G|=\alpha$ and
taking $K:=\Hom(G,\TT)=\widehat{G_d}$) that $w(G,\sT_H)=(\log\,\alpha)^\omega$
and $(G,\sT_H)$ has no convergent
sequences. (We note in this connection that
$(\log\,\alpha)^\omega<2^\alpha$ for every strong limit cardinal
$\alpha$ such that $\cf(\alpha)>\omega$.) This discussion
suggests a question.
}
\end{remark}

\begin{question}
{\rm
Given an infinite Abelian group $G$, what is the minimal
weight of a topology in $\sT\in\ttt(G)$ such that no nontrivial sequence
converges in $(G,\sT)$? For which $G$ is this $2^{|G|}$?
}
\end{question}


\begin{thebibliography}{99}

\bibitem{Barb} {\sc Barbieri, G., D. Dikranjan, C. Milan} and {\sc H. Weber,}
	{\em Answer to Raczkowski's questions on convergent sequences of integers,}
	Topology Appl. {\bf 132} (1) (2003), 89-101.
	
\bibitem{cantor84}
{\sc Georg Cantor}, \emph{{\"{U}}ber unendliche, lineare {P}unktmannigfaltigkeiten
  {V}{I}}, Math. Annalen \textbf{23} (1884), 453--458.


 \bibitem{WCSN:74} {\sc W. W. Comfort} and {\sc S. Negrepontis}, 
                  {\em The Theory of Ultrafilters,} 
                 Grundlehren der mathematischen Wissenschaften, vol. 211, Springer Verlag,
                  Berlin-Heidelberg-New York, 1974.

\bibitem{wclr:82} {\sc W. W. Comfort} and {\sc L. C. Robertson}, 
                  {\em Proper pseudocompact extensions of compact 
                  Abelian  group topologies,} 
                  Proc. Amer. Math. Soc., 
                  {\bf 86} (1) (1982), 173-178.
        
\bibitem{wclr:85} {\sc W. W. Comfort} and {\sc L. C. Robertson}, 
           {\em Cardinality constraints for pseudocompact and for 
           totally dense subgroups of compact topological groups,} 
           Pacific J. of Math., {\bf 119} (2) (1985), 265-285.
           
\bibitem{wclr:88} {\sc W. W. Comfort} and {\sc L. C. Robertson}, 
                  {\em Extremal phenomena in certain classes of totally bounded 
                  groups,} Dissertationes Math., PWN, {\bf 272} (1988).

\bibitem{wckr:64} {\sc W. W. Comfort} and {\sc K. A. Ross}, 
                  {\em Topologies induced by groups of characters,} 
                   Fundamenta Math., {\bf 55} (1964), 283-291.
                   MR {\bf 30}:183
                   
 \bibitem{wckr:66} {\sc W. W. Comfort} and {\sc K. A. Ross}, 
                  {\em Pseudocompactness and uniform continuity in topological groups,} 
                   Pacific J. Math., {\bf 16} (1966), 483-496.
 
                   
\bibitem{CTWu:93} {\sc W. W. Comfort, F. J. Trigos-Arrieta} and {\sc T.
	S. Wu}, {\em The Bohr compactification, modulo a metrizable subgroup},
	Fundamenta Math., {\bf 143} (1993), 119-136.
	Correction: same journal {\bf 152} (1997), 97-98. MR:94i22013,
	Zbl. 81222001.
	

\bibitem{diksha03}
{\sc Dikran~N. Dikranjan} and {\sc Dmitri~B. Shakhmatov}, \emph{Forcing hereditarily
  separable compact-like group topologies on {A}belian groups} (2003).
Manuscript submitted for publication.


\bibitem{vdi}
{\sc E.~K.~van Douwen}, \emph{The product of two countably compact topological
  groups}, Trans. Amer. Math. Soc. \textbf{262} (1980), 417--427.


\bibitem{vdie}
{\sc E.~K.~van Douwen}, \emph{The weight of a pseudocompact (homogeneous) space whose
  cardinality has countable cofinality}, Proc. Amer. Math. Soc. \textbf{80}
  (1980), 678--682.

\bibitem{efimov69}
{\sc B.~Efimov}, \emph{On imbedding of {S}tone-{\v{C}}ech compactifications of
  discrete spaces in bicompacta}, Soviet Math. Doklady \textbf{10} (1969),
  1391--1394. [Russian original in: {\cyr Doklady Akad. Nauk SSSR} {\bf 187}
  (1969), 244--266.]

\bibitem{engel}
{\sc Ryszard Engelking}, \emph{General {T}opology},
Heldermann Verlag, Berlin, 1989.

\bibitem{fedorcuk76}
{\sc V.~V. Fedor{\v{c}}uk}, \emph{Fully closed mappings and consistency of some
  theory of general topology with the ansioms of set theory}, Math.
  {U}{S}{S}{R} {S}bornik \textbf{28} (1996), 1--26. Russian original in:
{\cyr Matem. Sbornik} (N.S.) {\bf 99} (1976), 3--33.

\bibitem{fuch:70} {\sc L. Fuchs,} {\em Infinite Abelian Groups}, {\bf
	vol. I.} Academic Press. New York-San Francisco-London, 1970.




\bibitem{glickii}
{\sc Irving Glicksberg}, \emph{Uniform boundedness for groups}, Canadian J. Math.
  \textbf{14} (1962), 269--276.

\bibitem{Ha} {\sc P. Halmos}, {\em Measure Theory}, D. Van Nostrand, 
	New York, 1950.
	
\bibitem{HaKu} {\sc J. Hart} and {\sc K. Kunen}, {\em Limits in function
	spaces and compact groups,} Topology and Its Applications. To appear.

\bibitem{hartvmii}
{\sc K.~P. Hart} and {\sc J.~van Mill}, \emph{A countably compact group {$H$} such that
  {$H$}$\times${$H$} is not countably compact}, Trans. Amer. Math. Soc.
  \textbf{323} (1991), 811--821.
 
\bibitem{haus14}
{\sc Felix Hausdorff}, \emph{Grundz{\"u}ge der {M}engenlehre}, Veit, Leipzig, 1914,
  [Reprinted: Chelsea Publ. Co., New York, 1949.]

\bibitem{ehkr:63} {\sc E. Hewitt} and {\sc K. A. Ross}, {\em Abstract Harmonic 
         Analysis,} {\bf vol. I}. Springer Verlag, Berlin $\cdot$ G\"ottingen 
         $\cdot$ Heidelberg, 1963. MR {\bf 28}:58.


\bibitem{kakkod}
{\sc Shizuo Kakutani and Kunihiko Kodaira}, \emph{{\"{U}}ber das {H}aarsche {M}ass in
  der lokal bikompacten {G}ruppen}, Proc. Imperial Acad. Tokyo \textbf{20}
  (1944), 444--450. [Reprinted In: Selected papers of Shizuo Kakutani volume 1,
  edited by Robert R. Kallman, pp.~68--74. Birkh{\"{a}}user,
  Boston-Basel-Stuttgard, 1986.]

\bibitem{kodaira41}
{\sc Kunihiko Kodaira}, \emph{{\"{U}}ber die beziehung zwischen den {M}assen und den
  {T}opologien in einer {G}ruppe}, Proc. Physico-Math. Soc. Japan \textbf{16
  (Series 3)} (1941), 67--119.

\bibitem{VK:35} {\sc E. R. van Kampen}, {\em Locally bicompact Abelian groups
         and their character groups,} Annals of Math., {\bf 36} (2) (1935), 
         448-463.

\bibitem{leptin55a}
{\sc Horst Leptin}, \emph{Abelsche {G}ruppen mit kompakten {C}haraktergruppen und
  {D}ualit{\"{a}}tstheorie gewisser linear topologischer abelscher {G}ruppen},
  Abhandlungen Mathem. Seminar Univ. Hamburg \textbf{19} (1955), 244--263.

\bibitem{malyshap}
{\sc V.~I. Malykhin} and {\sc L.~B. Shapiro}, \emph{Pseudocompact groups without convergent
  sequences}, Mathematical Notes \textbf{37} (1985), 59--62. [Russian original
  in: {\cyr Matematicheskie Zametki {\bf 37} (1985), 103--109}.]

\bibitem{PvKd:34} {\sc L. Pontryagin}, {\em The theory of topological 
         commutative groups},  Annals of Math., {\bf 35} (2) (1934), 361-388.

\bibitem{raczphd}
{\sc Sophia~U. Raczkowski-Trigos,} \emph{Totally {B}ounded {G}roups}, Ph.D. thesis,
   Wesleyan University, Middletown, Connecticut, USA, 1998.

\bibitem{racz:02} {\sc S. U. Raczkowski}, {\em Totally bounded topological
	group topologies on the integers}, Topology Appl., 
	{\bf 121} (2002), 63-74. 
	
\bibitem{SaStr:85} {\sc S. Saeki} and {\sc K. R. Stromberg,} {\em Measurable subgroups
	and non-measurable characters,} Math. Scandinavica {\bf 57} (1985), 
	359-374.
	
\bibitem{shakh92}
{\sc Dmitri~B. Shakhmatov}, \emph{Compact spaces and their generalizations},
in: Recent
  {P}rogress in {G}eneral {T}opology (M.~Hu\v{s}ek and Jan van Mill, eds.),
pp.~571--640. North-Holland, Amsterdam-London-New York-Tokyo, 1992.

\bibitem{shak94}
{\sc Dmitri~B. Shakhmatov}, \emph{A direct proof that every infinite compact group
  {$G$} contains $\{0,1\}^{w(G)}$}, In: Papers on General Topology and
  Applications, {\it Annals of the New York Academy of Sciences}
vol.~{\bf 728} (Susan Andima, Gerald Itzkowitz, T.~Yung Kong, Ralph Kopperman,
  Prabud~Ram Misra, Lawrence Narici, and Aaron Todd, eds.),
pp.~276-283. New York, 1994. [Proc. June,
  1992 Queens College Summer Conference on General Topology and Applications.]

\bibitem{sirota}
{\sc S.~M. Sirota}, \emph{The product of topological groups and extremal
  disconnectedness}, Math. USSR Sbornik \textbf{8} (1969), 169--180. [Russian
  original in: {\cyr Matem. Sbornik {\bf 79} (121) (1969), 179--192}.]


		
\bibitem{steinhaus20}
{\sc Hugo Steinhaus}, \emph{Sur les distances des points des ensembles de measure
  positive}, Fund. Math. \textbf{1} (1920), 93--104.

\bibitem{strom72}
{\sc Karl~R. Stromberg}, \emph{An elementary proof of {S}teinhaus's theorem}, Proc.
  Amer. Math. Soc. \textbf{36} (1972), 308.
  
\bibitem{Str92} {\sc Karl~R. Stromberg,} {\em Universally nonmeasurable subgroups
	of $\RR$,} Math. Assoc. of Amer., {\bf 99} (3) (1992), 253-255.

\bibitem{tomita96}
{\sc Artur~Hideyuki Tomita}, \emph{On finite powers of countably compact groups},
  Comment. Math. Univ. Carolin. \textbf{37} (1996), 617--626.

\bibitem{tomita99}
{\sc Artur~Hideyuki Tomita}, \emph{A group under {\rm MA}$_{\rm countable}$ whose square is countably
  compact but whose cube is not}, Topology and Its Applications \textbf{91}
  (1999), 91--104.


\bibitem{Weil:37} {\sc Andr\'{e} Weil}, {\em Sur les {E}spaces \`{a} {S}tructure
                   {U}niforme et sur la {T}opologie {G}\'{e}n\'{e}rale}, Publ. Math. 
                   Univ. Strasburg, Hermann, Paris, 1937.
	
\bibitem{weilii}
{\sc Andr\'{e} Weil}, \emph{L'int\'{e}gration dans les {G}roupes {T}opologiques et
  ses {A}pplications}, Actualit\'{e}s Scientifiques et Industrielles \#869,
  Publ. Math. Institut Strasbourg, Hermann, Paris, 1940. [Deuxi\`eme \'edition
  \#1145, 1951.]

\end{thebibliography}
\end{document}